\documentclass[11pt]{amsart}

\usepackage[T1]{fontenc}
\usepackage{lmodern}
\usepackage[letterpaper, margin=1.25in]{geometry}
\usepackage{amsmath, amssymb, amsthm, mathtools}

\usepackage{enumitem}
\setlist[enumerate,1]{label=\textup{(\alph*)}, ref=(\alph*), leftmargin=*}

\usepackage[hidelinks]{hyperref}
\usepackage[capitalise, noabbrev, nameinlink]{cleveref}

\numberwithin{equation}{section}

\theoremstyle{plain}

\newtheorem{maintheorem}{Theorem}

\newtheorem{theorem}[equation]{Theorem}
\newtheorem{proposition}[equation]{Proposition}
\newtheorem{lemma}[equation]{Lemma}

\theoremstyle{definition}

\newtheorem{openproblem}{Open problem}

\theoremstyle{remark}

\crefname{maintheorem}{Theorem}{Theorems}

\newcommand{\R}{\mathbb R}
\newcommand{\C}{\mathbb C}

\DeclareMathOperator{\supp}{supp}

\DeclareMathOperator{\spec}{spec}
\renewcommand{\Re}{\operatorname{Re}}
\renewcommand{\Im}{\operatorname{Im}}

\usepackage{xcolor}
\colorlet{romain}{green!50!black}
\colorlet{josef}{blue!50!black}

\title{Inferring geometric structure from eigenfunction triple products}

\author[J. Greilhuber]{Josef Greilhuber}
\address{Department of Mathematics,
  Massachusetts Institute of Technology,
  77 Massachusetts Avenue,
  Cambridge, Massachusetts 02139,
  United States of America}
\email{jgreil@mit.edu}

\author[R. Speciel]{Romain Speciel}
\address{Department of Mathematics,
  Stanford University,
  450 Jane Stanford Way, Building 380,
  Stanford, California 94305,
  United States of America}
\email{rspeciel@stanford.edu}

\date{November 1, 2026}

\begin{document}

\begin{abstract}
  This article presents two converses to well-known observations regarding products of eigenfunctions of the Laplacian on a closed Riemannian manifold $(M,g)$. First, if $(M,g)$ is real-analytic, the coefficients in the eigenfunction expansion of a product of eigenfunctions decay exponentially. We show that such exponential decay in fact implies that $g$ is real-analytic. Second, it is clear that if $(M,g)$ is flat, there exists an eigenbasis for which the product of any two of its elements lies in the span of at most $N = N(M,g)$ elements. We prove that any manifold with such an eigenbasis must be flat. We then construct metrics on $\mathbb{S}^7$ that share an eigenbasis with the standard metric yet are not locally homogeneous, from which we deduce that homogeneity cannot be characterized solely in terms of the multiplicative structure of eigenfunctions.
\end{abstract}

\maketitle

\section{Introduction}
\label{sec:introduction}

We consider the Laplacian $\Delta$ on $(M,g)$, a smooth closed connected Riemannian manifold of dimension $n\geq 2$. Choose a sequence of Laplace eigenfunctions $(\phi_j)_{j=0}^\infty$ which forms an orthonormal basis of $L^2(M)=L^2(M;\mathbb{C})$ and is ordered such that, with $\Delta\phi_j=\mu_j^2\phi_j$, we have $0=\mu_0<\mu_1\leq \mu_2\leq \dots\to \infty$. The \emph{triple product constants} associated to this basis, defined by
\[
    c_{jk\ell}=\langle \phi_j\phi_k,\phi_\ell\rangle =\int_M \phi_j\phi_k\overline{\phi_\ell}\,d\textrm{vol}_g,
\]
encode the multiplicative structure of the eigenfunctions, as $\phi_j\phi_k=\sum_\ell c_{jk\ell}\phi_\ell$. While it is well known that the Laplace spectrum $(\mu_j^2)_{j=0}^\infty$ alone does not determine $(M,g)$ \cite{GordonPerrySchueth2005}, the triple product constants provide sufficient additional spectral information to allow for a complete reconstruction of the underlying manifold. Indeed, $(M,g)$ is uniquely recoverable from the data $(\mu_j, c_{jk\ell})$ in any orthonormal eigenbasis  \cite{Schaefer2024, SchildkrautSpeciel2026}.

In this paper, we refine this correspondence to better understand how eigenvalues and triple product constants encode the geometric features of the underlying manifold. To this end, we identify simple conditions on $(\mu_j, c_{jk\ell})$ which characterize real-analyticity or flatness of the metric. We then study counterexamples to a natural conjecture on how the structure of the triple product constants may determine locally homogeneous compact spaces.

\subsection{Results}
Our first two main results provide conditions on the spectrum and the triple product constants which identify exactly when a manifold is real-analytic or when it is flat. More precisely, calling a metric \emph{real-analytic} if the underlying manifold admits a real-analytic atlas, compatible with its smooth structure, in which the coefficients of the metric are real-analytic, we have the following.

\begin{maintheorem}\label{thm:realanalytic}
    Let $(M,g)$ be a smooth closed connected Riemannian manifold, and fix an orthonormal eigenbasis with corresponding eigenvalues $\mu_j^2$ and triple product constants $c_{jk\ell}$. The metric $g$ is real-analytic if and only if for each $j$, there are constants $a_j,C_j>0$ such that
    \[
        |c_{jk\ell}|\leq C_je^{-a_j|\mu_k-\mu_\ell|}
    \]
    for every $k,\ell$.
\end{maintheorem}
\noindent The main novel contribution here is the fact that this exponential decay of the triple product constants implies the underlying metric is real-analytic. \cref{thm:tech_real_analytic} provides a broader statement of this result along with two other equivalent characterizations, one in terms of how much a frequency exceeds its expected threshold, the other in terms of approximations by low-frequency eigenfunctions.

Similarly, we characterize flat manifolds purely in terms of their triple product constants, this time without reference to the eigenvalues. Recall that the \emph{holonomy group} of $M$ consists of those linear transformations on the tangent space at a fixed point which are induced by parallel transport along piecewise smooth loops based at that point. Up to conjugacy, the holonomy group does not depend on the base point in its definition.

\begin{maintheorem}\label{thm:flat}
    Let $(M,g)$ be a smooth closed connected Riemannian manifold. The metric $g$ is flat if and only if there exists an orthonormal eigenbasis, with corresponding triple product constants $c_{jk\ell}$, and an integer $N$ such that
    \[
        \#\{\ell:c_{jk\ell}\neq 0\}\leq N
    \]
    for every $j, k$. When $g$ is flat, the least such $N$, across all possible choices of orthonormal eigenbases, is the order of the holonomy group of $(M,g)$.
\end{maintheorem}
\noindent In other words, if there exists some $N$ such that the product of any two eigenbasis elements may be written as a linear combination of no more than $N$ basis elements, then $g$ is flat. This result is stated more broadly in \cref{thm:tech_flat}, which also allows us to count eigenspaces over which products split, rather than just eigenbasis terms.

The key initial step in the proofs of both of these results is an application of Egorov's theorem which converts the structure of triple product constants into control of eigenfunctions along geodesics on the manifold. With $M_f$ denoting the multiplication operator by some function $f\in C^\infty(M)$, we study the operator
\begin{equation}\label{eq:conjugated_multiplication}
    A_f(t)=e^{it\sqrt{\Delta}}\,M_f\,e^{-it\sqrt{\Delta}}.
\end{equation}
Egorov's theorem relates the matrix entries of $A_f(t)$ to the composition of $M_f$ with the geodesic flow. We apply this fact to deduce \cref{prop:decay-to-geodesic} and \cref{prop:sparse-to-geodesic}, which underpin the recovery results.

In the real-analytic case, we obtain from this analytic coordinates, from which the structure of the metric follows. Conversely, when the metric is real-analytic, we show eigenfunctions extend holomorphically to a Grauert tube, then recover the decay condition. In the flat case, we construct a torus action on an open subset of the unit sphere bundle, from which we deduce flatness after applying tools from semiclassical analysis.

Finally, we consider triple product constants on compact homogeneous spaces. Denote the algebraic span of the eigenfunctions on $M$ by $\mathcal{A}(M,g)$. In other words, $\mathcal{A}(M,g)$ consists of all finite sums of eigenfunctions on $M$. A simple application of the Peter-Weyl theorem shows that, on compact homogeneous spaces, the product of any two eigenfunctions is always a finite sum of eigenfunctions \cite[Chapter III, Sections 1 and 3]{BrockerTomDieck1985}. Equivalently, when $M$ is compact and homogeneous, $\mathcal{A}(M,g)$ forms an algebra. Note this proprety passes to finite covers. It is therefore natural to ask, in analogy to \cref{thm:realanalytic} and \cref{thm:flat}, if the property that $\mathcal{A}(M,g)$ is an algebra characterizes those spaces. However, this is false. In fact, the following theorem shows the situation to be quite dire.

\begin{maintheorem}\label{thm:nonhomogeneous}
    On $\mathbb{S}^7$, there exists a metric which is not locally homogeneous and which shares an orthonormal eigenbasis with the standard metric. In particular, triple product constants alone cannot detect local homogeneity.
\end{maintheorem}

\noindent This result is stated in more generality in \cref{thm:tech_nonhomogeneous}, which is proved by rescaling the standard metric along three orthogonal Killing fields then showing that a resulting eigenspace must be one dimensional, thus ruling out homogeneity.

\subsection{Open problems} We present here several natural directions for further investigation which refine the relationship between the geometry of a manifold and the structure of its eigenvalues and triple product constants.

\begin{openproblem}\label[openproblem]{prob:finite_generation}
    \emph{Prove or disprove: if $\mathcal{A}(M,g)$ is an algebra, then it is finitely generated.} A related finite-generation question for Laplacian subalgebras on compact homogeneous spaces appears in \cite[Question 2]{LinMendesRadeschi2025}. On compact homogeneous spaces, one can show that the algebra generated by eigenfunctions is indeed finitely generated. However, the methods employed to prove this do not seem to extend to the general case.
\end{openproblem}

\begin{openproblem}
    \emph{Let $M$ be a surface. Prove or disprove: $\mathcal{A}(M,g)$ is an algebra if and only if $(M,g)$ has constant nonnegative curvature.} One direction is immediate. Resolving this question would also resolve \cref{prob:finite_generation} in the case of surfaces.
\end{openproblem}

\begin{openproblem}
    \emph{Prove or disprove: If $\mathcal{A}(M,g)$ is an algebra, then the geodesic flow on $S^\ast M$ is not ergodic.} One may hope to find an application of Egorov's theorem as in the proofs of \cref{thm:realanalytic,thm:flat} to display eigenfunctions whose product spreads over the entire spectrum.
\end{openproblem}

\subsection{Structure}
In \cref{sec:background}, the main theorems are restated in greater generality and related works are discussed. The remainder of the paper develops the proofs of these statements. \cref{sec:egorov} applies Egorov's theorem to convert conditions on the triple product constants into properties of eigenfunctions along geodesics. This is employed to characterize real-analyticity in \cref{sec:real-analyticity} and flatness in \cref{sec:flatness}. Finally, in \cref{sec:weird-spheres}, we construct metrics on $\mathbb{S}^7$ which are not locally homogeneous yet whose triple product constants exhibit a certain rigid structure.

\subsection{Acknowledgments} The authors thank Anshul Adve, Ben Church, Elena Kim, Eugenia Malinnikova, Rafe Mazzeo, Zeev Rudnick, Carl Schildkraut, and Mitchell Taylor for many enlightening conversations and much wise advice.

This project was completed in collaboration with several AI models, primarily including ChatGPT 5.6 Sol and ChatGPT 6 Astra. The authors first considered the questions addressed in this paper in early 2025, then made incremental progress, mostly without AI assistance, until the spring of 2026. This progress included a partial result towards \Cref{thm:flat}, as described in \cite{SchildkrautSpeciel2026}, and \Cref{prop:sparse-to-geodesic}. At this point, through extensive mathematical conversations with AI models, the authors were able to obtain the proofs of the main results. In particular, the models suggested the proof of \Cref{prop:decay-to-geodesic}, the main ideas behind the arguments of \Cref{sec:flatness}, and the family of counterexamples studied in \Cref{thm:tech_nonhomogeneous}. The proofs suggested by the models were correct, but relied on often convoluted auxiliary constructions and usually did not impart geometric intuition to their readers. The authors worked to significantly simplify these proofs, in the process finding various alternative arguments. Then, they dedicated significant efforts to finding clear and elegant ways of presenting the results, and accurately attributing the employed ideas. Towards the end of the writing process, the authors also employed AI tools (both ChatGPT and Claude) for minor editorial purposes. This manuscript is human written, in its entirety, by the authors, who bear all responsibility for its content and correctness.

\section{Background and technical statements}\label{sec:background}

Eigenfunction triple products have long been studied in number theory due to their connection to $L$-functions, as expressed by Watson's formula \cite[Theorem 3]{Watson2002}. In this setting, exponentially decaying upper bounds on the triple product constants were established by Sarnak \cite{Sarnak1994} and Bernstein and Reznikov \cite{BernsteinReznikov1999}, while complementary averaged lower bounds were established by Reznikov \cite[Theorem 2]{Reznikov2001}. More recently, the triple product constants have appeared in work of Adve on the conformal bootstrap \cite{Adve2025}, and in work of Bisain, Humphries, Mandelshtam, Walsh, and Wang on quantum ergodicity and quantum chaos \cite{BisainHumphriesMandelshtamWalshWang2024}.

We employ triple product constants to characterize real-analyticity through the following theorem, which generalizes \cref{thm:realanalytic}.

\begin{theorem}
\label[theorem]{thm:tech_real_analytic}
    Let $(M,g)$ be a smooth closed connected Riemannian manifold, and fix an orthonormal eigenbasis with corresponding eigenvalues $\mu_j^2$ and triple product constants $c_{jk\ell}$. The following are equivalent:

    \begin{enumerate}

        \item The metric $g$ is real-analytic.

         \item \label{cond:analytic-triangle} There are constants $a,C>0$, $b\geq 0$ such that, with
        \[
            R_{jk\ell}=\max\{0,\mu_j- (\mu_k+\mu_\ell),\, \mu_k-(\mu_\ell+\mu_j),\, \mu_\ell-(\mu_j+\mu_k)\},
        \]
        we have
        \[
             |c_{jk\ell}|\leq
             C(1+\mu_j+\mu_k+\mu_\ell)^b
             e^{-a R_{jk\ell}}
        \]
        for every $j,k,\ell$.

        \item \label{cond:analytic-tail} There are constants $a,C>0$, $b\geq0$ such that, for all $R\geq0$, we have
        \[
            \Bigl(\sum_{\substack{\ell\textup{ with}\\ \mu_\ell>
            \mu_j+\mu_k+R}}|c_{jk\ell}|^2\Bigr)^{1/2}
            \leq C(1+\mu_j+\mu_k)^b e^{-aR}
        \]
        for every $j,k$.  
        
        \item \label{cond:analytic-gap} For each $j$, there are constants $a_j,C_j>0$, $b_j\geq 0$ such that
        \[
            |c_{jk\ell}|\leq C_j(1+\mu_j+\mu_k)^{b_j}e^{-a_j|\mu_k-\mu_\ell|}
        \]
        for every $k,\ell$.

    \end{enumerate}

    \noindent When $g$ is real-analytic, one may choose $a_j$ independently of $j$ and $b_j=0$ in condition \ref{cond:analytic-gap}.
    
\end{theorem}

\noindent We shall prove this result by showing (a) $\implies$ (b) $\implies$ (c) $\implies$ (d) $\implies$ (a). Condition \ref{cond:analytic-triangle} expresses the decay of the triple product constants in terms of how much a frequency exceeds its expected threshold, as suggested by the intuition that frequency should add when eigenfunctions are multiplied \cite{Steinerberger2019}. This perspective is explored in detail in work of Wyman \cite{Wyman2022}. Condition \ref{cond:analytic-tail} may be understood as a statement about how accurately the products of eigenfunctions are approximated by low dimensional spaces, as studied by Lu, Sogge, and Steinerberger in \cite[Theorem 1]{LuSoggeSteinerberger2019}. Finally, condition \ref{cond:analytic-gap} is analogous to Seeley's characterization of analytic functions in \cite{Seeley1969}. Zelditch \cite[Proposition 1]{Zelditch2012} establishes exponential decay for fixed input eigenfunctions, while Charron and Pagano \cite[Theorem 1.1]{CharronPagano2025} give quantitative exponential bounds beyond a manifold-dependent multiple of the input frequencies. Our main contribution is the converse implication, namely that one can recover real-analyticity from the estimates in this theorem.

Next, we address the case of flat manifolds by proving the following theorem, which generalizes \cref{thm:flat}.

\begin{theorem}
\label[theorem]{thm:tech_flat}
    Let $(M,g)$ be a smooth closed connected Riemannian manifold. The following are equivalent:
    \begin{enumerate}
        \item The metric $g$ is flat.
        \item\label{cond:bounded-terms} There exists an orthonormal eigenbasis, with corresponding triple product constants $c_{jk\ell}$, and an integer $N\geq 1$ such that
        \[
            \#\{\ell:c_{jk\ell}\neq 0\}\leq N
        \]
        for every $j, k$.
        \item \label{cond:bounded-espace} There exists an orthonormal eigenbasis, with corresponding triple product constants $c_{jk\ell}$, and an integer $N\geq 1$ such that
        \[
            \#\{\mu_\ell:c_{jk\ell}\neq 0\}\leq N
        \]
        for every $j, k$.
    \end{enumerate}
    When $g$ is flat, the least $N$ across all possible choices of orthonormal eigenbases for which condition \ref{cond:bounded-terms} holds is the same as that for which condition \ref{cond:bounded-espace} holds. This number is the order of the holonomy group of $(M,g)$.
\end{theorem}

\noindent Note that condition \ref{cond:bounded-espace} counts only the number of eigenspaces which meet a particular product, rather than the number of basis terms needed to write this product. This result, and a generalization to the graph setting, was proven in the case of real-valued eigenfunctions when $N=2$ by Schildkraut and the second author in \cite{SchildkrautSpeciel2026}. The theorem above allows for any uniform bound on the number of terms.

We conclude by proving the following theorem, which strengthens \cref{thm:nonhomogeneous}.

\begin{theorem}
\label[theorem]{thm:tech_nonhomogeneous}
    There exists a family of real-analytic metrics $g_t$ on $\mathbb{S}^7$, depending real-analytically on the parameter $t$ defined in a neighborhood of $1$, such that:
    \begin{enumerate}
        \item the metric $g_1$ is the standard metric;
        \item these metrics induce the same volume form;
        \item these metrics share a common orthonormal eigenbasis;
        \item and the metrics $g_t$ are not locally homogeneous when $t\neq 1$.
    \end{enumerate}
    In particular, the triple product constants corresponding to the basis obtained by condition (c) are identical across this family.
\end{theorem}

\noindent Since these metrics all share a common orthonormal eigenbasis, we have that, for $t$ sufficiently close to $1$, the ring $\mathcal{A}(M,g_t)$ consists of the restrictions of polynomials in $\R^8$ to the sphere. Furthermore, since $g_t\to g_1$ in $C^\infty$ as $t\to 1$, one cannot hope for any local stability statement purely in terms of the triple product constants. This theorem sheds light on the question whin appears in \cite[Question 3]{LinMendesRadeschi2025} of understanding those closed manifolds for which $\mathcal{A}(M,g)$ is an algebra. In particular, such manifolds need not be locally homogeneous.

\section{Eigenfunctions along geodesics}\label{sec:egorov}

In this section, we relate the structure of the triple product constants to the behavior of eigenfunctions along geodesics. Our first proposition addresses the case in which the triple product constants decay exponentially.

\begin{proposition}\label[proposition]{prop:decay-to-geodesic}
    Fix a function $f\in C^\infty(M)$. If there exists an orthonormal eigenbasis $(\phi_j)_{j=0}^\infty$ and constants $a, C>0$, $b\geq 0$ such that 
    \[
        |\langle f\phi_j, \phi_k\rangle |\leq C(1+\mu_j+\mu_k)^b e^{-a|\mu_j-\mu_k|}
    \]
    for every $j,k$, then there exists some $r>0$ such that, for every unit-speed geodesic $\gamma$, the function $f\circ \gamma$ extends holomorphically to the strip $\{z\in \mathbb{C}\,\colon\,|\Im(z)|<r\}$, on which it is bounded independently of $\gamma$.
\end{proposition}

\noindent Our second proposition treats the case in which the triple product constants are sparse.  

\begin{proposition}\label[proposition]{prop:sparse-to-geodesic}
    Fix a function $f\in C^\infty(M)$. If there exists an orthonormal eigenbasis $(\phi_j)_{j=0}^\infty$ and a constant $N>0$ such that 
    \[
        \#\{\mu_k\,:\, \langle f\phi_j,\phi_k\rangle\neq 0 \}\leq N,
    \]
    for every $j$, then, for every geodesic $\gamma$, the function $f\circ \gamma$ is a linear combination of no more than $N$ complex exponentials.
\end{proposition}

Both corresponding proofs employ Egorov's theorem to translate the multiplicative properties of the eigenfunction into dynamical properties on the underlying manifold. We therefore begin by recalling its statement, in the simplified form which we shall employ, along with two related facts. We then proceed to prove our propositions.

Denote by $\Psi^0_\textrm{cl}(M)$ the space of classical pseudodifferential operators of order zero on $M$, and write $\sigma(A)$ for the principal symbol of $A\in\Psi^0_\textrm{cl}(M)$. For any $t\in \R$, set $G_t\colon S^*M\to S^*M$ to be the time-$t$ geodesic flow.

\begin{theorem}[Egorov's Theorem {\cite[Theorem 9.1]{Zelditch2017}}]
    For $A\in \Psi^0_\textrm{cl}(M)$ and $t\in \R$, we have
    \[
        A(t)=e^{it\sqrt{\Delta}}Ae^{-it\sqrt{\Delta}}\in \Psi^0_\textrm{cl}(M)\quad \textrm{and}\quad \sigma(A(t))=\sigma(A)\circ G_t.
    \]
\end{theorem}

\noindent In particular, with $A_f(t)$ as in \Cref{eq:conjugated_multiplication} and $\pi\colon S^*M\to M$ the standard projection, the theorem gives
\begin{equation}
    \sigma(A_f(t))=f_t,\quad\textrm{where}\quad f_t=f\circ \pi \circ G_t.
\end{equation}

We also recall that any sequence $(\phi_j)_{j=0}^\infty$ of normalized eigenfunctions with $\mu_j\to \infty$ admits a subsequence $(\phi_{j_k})_{k=0}^\infty$ and a corresponding probability measure $\mu$ on $S^*M$ such that, for any $A\in \Psi^0_\textrm{cl}(M)$,
\begin{equation}
    \langle A\phi_{j_k},\phi_{j_k} \rangle\to \int_{S^*M}\sigma(A)\,d\mu
\end{equation}
as $k\to \infty$ \cite[Section 5.1]{Zworski2012}. We call $\mu$ the \emph{semiclassical measure} corresponding to the subsequence $(\phi_{j_k})_{k=0}^\infty$. A simple application of Egorov's Theorem, after noting that 
\[
    \langle A\phi_{j_k},\phi_{j_k} \rangle=\langle e^{it\sqrt{\Delta}}Ae^{-it\sqrt{\Delta}}\phi_{j_k},\phi_{j_k} \rangle,
\]
shows that all such measures must be $G_t$-invariant. Furthermore, semiclassical measures have the following property.

\begin{lemma}\label[lemma]{lem:semiclassical-measures-are-dense}
    Fix an orthonormal eigenbasis $(\phi_j)_{j=0}^\infty$. The union of the supports of the semiclassical measures corresponding to subsequences of $(\phi_j)_{j=0}^\infty$ is dense in $S^*M$.
\end{lemma}
\begin{proof}
    Pick $U\subset S^*M$ nonempty and open, and choose $A\in \Psi^0_\textrm{cl}(M)$ whose symbol is nonzero and supported in $U$. The local Weyl law states that
    \[
        \frac{1}{j}\sum_{k\leq j}\Vert{A\phi_k}\Vert^2\to L=\int_{S^*M} |\sigma(A)|^2\,d\mu_L>0
    \]
    as $j\to \infty$, where $d\mu_L$ is the standard normalized Liouville measure on $S^*M$ \cite[Theorem 8.7]{Zelditch2017}. It follows that $\Vert{A\phi_k}\Vert^2\geq L/2$ for infinitely many values of $k$. Any semiclassical measure obtained from a subsequence thereof must have support intersecting $U$.
\end{proof}

\noindent Armed with these facts, we may now proceed to the proofs of our two propositions.

\subsection{Decay of the triple product constants}
To establish \Cref{prop:decay-to-geodesic}, we begin with a fact from complex analysis.

\begin{lemma}\label[lemma]{lem:holomorphic_extention}
    Let $g\colon(-1,1)\to \mathbb{C}$ and suppose that, for some $N\geq 0$, there exist holomorphic functions $g_k$, $k\geq 1$, on $\mathbb{D}$ such that
    \[
        \sup_{z\in\mathbb{D}}|g_k(z)|\leq 2^{kN}\quad\textrm{and} \quad \sup_{t\in (-1,1)}|g_k(t)-g(t)|\leq 2^{-k}.
    \]
    Then, there exist constants $0<r<1$ and $M>0$, depending only on $N$, and a holomorphic function $h\colon\{|z|<r\}\to \mathbb{C}$ such that
    \[
        |h(z)|<M\quad\textrm{for}\quad |z|<r\qquad\textrm{and}\qquad  h(t)=g(t)\quad\textrm{for}\quad t\in (-r,r).
    \]
\end{lemma}

\begin{proof}
    Applying the assumptions to $h_k=g_k-g_{k-1}$, $k\geq 2$, gives
    \[
        \sup_{z\in\mathbb{D}}|h_k(z)|\leq 2\cdot 2^{kN}\quad\textrm{and} \quad \sup_{t\in (-1,1)}|h_k(t)|\leq 3\cdot 2^{-k}.
    \]
    By the maximum modulus principle, we have that for $0<\rho<1$ and $m= \lceil 2(N+1) \rceil $, the functions
    \[
        \frac{h_k(\rho z)}{(1-iz)^{km}} \quad\textrm{on the upper half disc}\qquad \textrm{and}\qquad\frac{h_k(\rho z)}{(1+iz)^{km}} \quad\textrm{on the lower half disc}
    \]
    are bounded by $3\cdot 2^{-k}$, since the denominators have norm at least $1$ on the real axis and at least $2^{km/2}$ on their respective semicircles. Taking $\rho\to 1$, we deduce
    \begin{equation}\label{eq:partial_sum_estimate}
        |h_k(z)|\leq 3\cdot 2^{-k}(1+|z|)^{km}
    \end{equation}
    for $z\in \mathbb{D}$.

    Pick $0<r<1$ sufficiently small so that $(1+r)^m/2<1$, noting this choice depends on $N$ only. The series
    \[
        h= g_1+\sum_{k=2}^\infty h_k
    \]
    converges uniformly on $\{|z|<r\}$, hence $h$ is holomorphic on $\{|z|<r\}$. The function $h$ agrees with $g$ on $t\in(-r,r)$ since the partial sums above give $g_k$, which agrees with $g$ arbitrarily well as $k\to \infty$, by assumption. The uniform bound on $h$ follows by summing the estimate of \Cref{eq:partial_sum_estimate}. This bound also depends on $N$ only.
\end{proof}

The following standard lemma provides a quantitative refinement of Egorov's Theorem, which we shall need in order to apply the previous lemma. We make use of the Sobolev spaces $H^s(M)$, $s\geq 0$, with norm
\[
    \lVert u\rVert_s^2=\lVert (1+\sqrt\Delta)^{s}u\rVert^2=\sum_{j=0}^\infty(1+\mu_j)^{2s}|\langle u, \phi_j\rangle|^2.
\]

\begin{lemma}\label[lemma]{lem:quantitative_egorov}
    Fix $(x_0,\xi_0)\in S^*M$ and choose real valued $\phi, \chi\in C^\infty(M)$ such that $\phi(x_0)=0$, $d\phi(x_0)=\xi_0$, and $\chi\equiv 1$ near $x_0$, with compact support on which $d\phi$ is nonvanishing. For $\kappa\geq 1$, define
    \[
        u_\kappa(x)=\chi(x)e^{i\kappa\phi(x)}.
    \]
    Then, for every $A\in \Psi^0_\textrm{cl}(M)$ and $s>0$, there exist constants $C_1$ and $C_2$ such that
    \[
        \lVert u_\kappa\rVert_s\leq C_1\kappa^s\quad\textrm{and}\quad \sup_{|t|\leq 1}|(A(t)u_\kappa)(x_0)-\sigma(A(t))(x_0,\xi_0)|\leq \frac{C_2}{\kappa}.
    \]
    The functions $\phi$ and $\chi$ may be chosen so that $C_1$ and $C_2$ do not depend on $(x_0,\xi_0)$.
\end{lemma}

\noindent This follows from the standard oscillatory testing formula for the principal symbol and Egorov's theorem, see \cite[Proposition 2.5 and Theorem 9.1]{Zelditch2017}.

We may now proceed to the proof of the proposition.

\begin{proof}[Proof of \Cref{prop:decay-to-geodesic}]
    When $u\in L^2(M)$ and $t\in \mathbb{R}$, we have
    \begin{equation}\label{eq:Af_written_out}
        A_f(t)u=\sum_{j,k}e^{it(\mu_k-\mu_j)}\langle u, \phi_j\rangle\langle f\phi_j,\phi_k\rangle \phi_k.
    \end{equation}
    We seek to control the right hand side of this equation so as to define a holomorphic extension of $A_f$. To achieve this, begin by noting that the operator $P_{jk}(u)= \langle u, \phi_j\rangle \phi_k$, when viewed as a map from $H^{3s}(M)\to H^s(M)$, has norm
    \[
        (1+\mu_j)^{-3s}(1+\mu_k)^s=\left(\frac{1+\mu_k}{1+\mu_j}\right)^{2s}\frac{1}{(1+\mu_j)^s(1+\mu_k)^s}.
    \]
    For $z\in \mathbb{C}$ with $|\Im(z)|<a/2$, we use
    \[
        \frac{1+\mu_k}{1+\mu_j}\leq 1+|\mu_k-\mu_j|\quad\textrm{and}\quad 1+\mu_j+\mu_k\leq (1+\mu_j)(1+\mu_k),
    \]
   along with our decay assumption for $\langle f\phi_j,\phi_k\rangle$, to find
    \begin{align*}
        \left| e^{iz(\mu_k-\mu_j)}\langle f\phi_j,\phi_k\rangle \lVert P_{jk}\rVert_{H^{3s}\to H^s} \right|&\leq  C e^{-a|\mu_k-\mu_j|/2}\left(\frac{1+\mu_k}{1+\mu_j}\right)^{2s}\frac{(1+\mu_j+\mu_k)^b}{(1+\mu_j)^s(1+\mu_k)^s}\\
        &\leq C' (1+\mu_j)^{b-s}(1+\mu_k)^{b-s}
    \end{align*}
    for some new constant $C'>0$. Therefore, for $s>b+n$, we have by Weyl's law that the sum in \Cref{eq:Af_written_out} converges uniformly as an operator from $H^{3s}(M)$ to $H^s(M)$ when taking $t=z$ with $|\Im(z)|<a/2$. Deduce $A_f(z)\colon H^{3s}(M)\to H^s(M) $ is holomorphic in $z$ on this strip, with uniformly bounded norm.

    Now, fix $(x_0,\xi_0)\in S^*M$, take $u_\kappa$ as in \Cref{lem:quantitative_egorov}, then define $g_\kappa(z)=(A_f(z)u_\kappa)(x_0)$ and $g(t)=\sigma(A_f(t))(x_0,\xi_0)$. Since we chose $s>b+n>n/2$, evaluation at a point is continuous on $H^s$ so $g_\kappa(z)$ is well defined when $|z|<a/2$. By \Cref{lem:quantitative_egorov}, we therefore obtain
    \[
        |g_\kappa(z)|\leq C_1'\kappa^{3s}\quad\textrm{and}\quad\sup_{|t|<1}|g_\kappa(t)-g(t)|<\frac{C_2'}{\kappa}.
    \]
    After rescaling our variables and applying \Cref{lem:holomorphic_extention}, we obtain constants $r,M>0$, independent of $(x_0,\xi_0)$, for which $g$ extends to $\{|z|<r\}$ on which it is bounded by $M$. The result follows after repeating this argument about $G_t(x_0,\xi_0)$ instead and applying Egorov's Theorem.
\end{proof}

\subsection{Sparsity of the triple product constants} We turn to our second proposition.

\begin{proof}[Proof of \Cref{prop:sparse-to-geodesic}]
    Let $P_k$ denote the projection onto the $\mu_k$-eigenspace of $\sqrt{\Delta}$, and observe that, under the assumptions in our statement, the sum
    \[
        A_f(t)\phi_j=\sum_{\mu_k} e^{it(\mu_k-\mu_j)}P_k(f\phi_j)
    \]
    has at most $N$ nonzero terms, by assumption. Therefore, for any choice $t_0,\dots, t_N\in \R$, there exists a linear dependence among the functions $A_f(t_k)\phi_j$, with $0\leq k\leq N$. It follows that the matrix
    \[
        \bigg(\langle A_f(t_k)\phi_j,A_f(t_\ell)\phi_j\rangle \bigg)_{k,\ell=0}^N
    \]
    has determinant $0$. Deduce, by taking limits along subsequences and applying Egorov's theorem, that, for any semiclassical measure $\mu$ arising from a subsequence of $(\phi_j)_{j=0}^\infty$, the matrix
    \[
        \Bigg(\int_{S^*M}f_{t_k}\overline{f_{t_\ell}}\,d\mu\Bigg)_{k,\ell=0}^N
    \]
    also has determinant zero.
    
    We deduce that there is therefore a nontrivial linear combination of the functions $(f_{t_k})_{k=0}^N$ which is identically zero in $L^2(S^*M,d\mu)$. Since $\mu$ is $G_t$-invariant, this same linear combination of the functions $(f_{s+t_k})_{k=0}^N$ must vanish in $L^2(S^*M,d\mu)$ for any $s\in \mathbb{R}$. This implies that, for any choice $s_0,\dots, s_N\in \R$ and $(x,\xi)\in \textrm{supp }\mu$, the determinant of the matrix
    \[
        \big(f_{s_k+t_\ell}(x,\xi)\big)_{k,\ell=0}^N
    \]
    is, once again, zero. This determinant is continuous in $(x,\xi)$ and the argument above holds for any semiclassical measure $\mu$. It follows by \Cref{lem:semiclassical-measures-are-dense} that the determinant vanishes identically on $S^*M$, for every choice of $s_k,t_\ell\in \R$. Conclude that for any fixed $(x,\xi)$, the space $\textrm{span}\{s\mapsto f_{s+t}(x,\xi)\,\colon\, t\in \mathbb{R}\}$ is of dimension at most $N$. A simple argument shows that any bounded continuous function on $\mathbb{R}$ whose translates span a space of dimension at most $N$ must be a linear combination of no more than $N$ complex exponentials, and the result follows.
\end{proof}

\section{Characterizing real-analyticity}\label{sec:real-analyticity}

We complete the proof of \Cref{thm:tech_real_analytic} in this section, by first obtaining a decay condition on the triple product constants from real-analyticity, then establishing implications between the several decay conditions presented in the theorem, and finally recovering real-analyticity from this decay.

\subsection{Decay from real-analyticity} To show that real-analyticity implies condition (b) in \Cref{thm:tech_real_analytic}, i.e., exponential decay of the triple product constants $\int_{M} \phi_j \phi_k \overline{\phi_\ell}\,d\mathrm{vol}_g$ away from the classically allowed region, it suffices to essentially repeat an argument from Zelditch's seminal article on pluripotential theory of Grauert tubes (\cite{Zelditch2012}, proof of Proposition 1).

%As shown by Whitney and Bruhat,
For the convenience of the reader, we repeat the necessary part of the exposition in Zelditch's article on the properties of the Grauert tube, referring the interested reader to \cite[Section 1]{Zelditch2012} for the infinitely better original. An $n$-dimensional real-analytic manifold $M$ embeds into a complex manifold $M_{\C}$, its so-called \emph{complexification}, as a maximally totally real submanifold. A real-analytic metric $g$ on $M$ gives rise to a distance function $r: M \times M \to \mathbb R$ such that $r^2$ is real-analytic near the diagonal $D = \{ (x,x) : x \in M \}$. By complexification, $r^2$ extends to a holomorphic function on a neighborhood of $D$ in $M_\C \times M_\C$, and by symmetry, $r^2(z,w) = r^2(w,z)$ for all $(z,w) \in M_\C \times M_\C$ where $r^2$ is defined. It follows that the function $\sqrt{\rho}$ defined by
\begin{align*}
    \sqrt{\rho} := \frac1{2i} \sqrt{r^2(z,\bar z)}
\end{align*}
is nonnegative near $M$ in $M_\C$ and vanishes precisely on $M$. The sublevel set
\begin{align*}
    M_\tau := \{z \in M_\C: \sqrt{\rho}(z) < \tau \}
\end{align*}
is a pseudoconvex domain in $M_\C$ called the \emph{Grauert tube}. For small $\tau$, its boundary $\partial M_\tau$ is diffeomorphic to the cosphere bundle $S^\ast M$. It turns out that, for $\tau$ smaller than a certain critical threshold $\tau_{an} > 0$, all eigenfunctions extend to $M_\tau$ simultaneously and satisfy the following uniform estimates.
\begin{theorem}[\cite{Zelditch2012}, Corollary 3, abridged]
    \label[theorem]{thm:Zelditch_cor_3}
    There exists $\tau_{an} > 0$, depending only on $M$ and $g$, such that every eigenfunction $\phi_j$ with $\Delta \phi_j = \mu_j^2 \phi_j$ extends analytically to $M_{\tau_{an}}$, with
    \begin{align*}
        \sup_{z \in M_\tau} |\phi_j(z)| \leq C_\tau \mu_j^{\frac{n+1}{4}} e^{\tau\mu_j}
    \end{align*}
    for every $\tau \leq \tau_{an}$.
\end{theorem}
We write $\mathcal A(\tau)$ for the operator of analytic continuation from $M$ to $\partial M_\tau$.

Consider the half-wave propagator $e^{i t \sqrt{\Delta}}: L^2(M) \to L^2(M)$, spectrally defined via
\begin{align*}
    e^{i t \sqrt{\Delta}} \phi_j = e^{i t \mu_j} \phi_j.
\end{align*}Its Schwartz kernel $U(t,\cdot,\cdot): \mathbb R \to \mathcal D'(M \times M)$ admits an analytic continuation to $U(t+i\tau,\zeta,x)$, where $\tau \leq \tau_{an}$, $\sqrt{\rho}(\zeta) \leq \tau_{an}$. The principal result on this operator, due to Boutet de Monvel \cite{BoutetdeMonvel1978}, is that it is \emph{invertible} even for non-real $t+i \tau$, provided one considers its image on the correct subset of $M_\C$ -- namely, the boundary of the Grauert tube. We write $\mathcal O^{s}(\partial M_\tau)$ for the intersection of the space $\mathcal O(\partial M_\tau)$ of distributional boundary values of holomorphic functions in $M_\tau$ with the Sobolev space $H^{s}(\partial M_\tau)$.
\begin{theorem}[\cite{Zelditch2012}, Theorem 7.1, abridged]
    For every $s$, the operator $U(i \tau)$ is a continuous isomorphism $H^s(M) \to \mathcal O^{s + \frac{n-1}{4}}(\partial M_\tau)$.
\end{theorem}

We are now prepared to expand upon Zelditch's argument to obtain exponential decay of the triple product constants $c_{jk\ell}$.

\begin{proof}[Proof of \Cref{thm:tech_real_analytic}] (a) $\Rightarrow$ (b).
Let $M_{\phi_j}: H^{s+\frac{n-1}{4}}(\partial M_\tau) \to H^{s+\frac{n-1}{4}}(\partial M_\tau)$ be the operator of multiplication by $\phi_j$. By the product rule, for a fixed sufficiently large integer $q$, 
$\|M_{\phi_j}\| \lesssim \|\phi_j^\C\|_{C^q(\partial M_\tau)}$. We combine Cauchy's estimates, applied to a slightly larger Grauert tube with $\tau' \approx \tau + \frac{1}{\mu_j}$, with the bounds of \Cref{thm:Zelditch_cor_3} to obtain
\begin{align*}
    \|\phi_j^\C\|_{C^q(\partial M_\tau)} \lesssim \mu_j^q \|\phi_j^\C\|_{C^0(\partial M_{\tau'})} \lesssim \mu_j^{q+\frac{n+1}{4}} e^{\tau' \mu_j} \lesssim \mu_j^{q+\frac{n+1}{4}} e^{\tau \mu_j}.
\end{align*}
Suppose $f$ is a function on $M$ which admits an analytic continuation to a neighborhood of $\overline{M_{\tau'}}$ for some $\tau' > 0$. Then the function 
\begin{align*}
    f_{t + i \tau} = \begin{cases} 
        U(t)^{-1} U(i \tau)^{-1} \mathcal A(\tau) f & 0 \leq \tau < \tau' \\ 
        e^{-i(t+i\tau) \sqrt{\Delta}} f & \tau < 0 \\ 
        \end{cases}
\end{align*}
provides an analytic continuation of $e^{- i t \sqrt{\Delta}} f$ to the half-plane $\{t+i \tau: \tau < \tau ' \}$. In particular, 
\begin{align*}
    \left \langle U(i \tau)^{-1} \mathcal A(\tau) f , \phi_j \right \rangle = e^{\tau \mu_j} \left \langle f , \phi_j \right \rangle.
\end{align*}
We therefore obtain 
\begin{align*}
    \mu_j^{q + \frac{n+1}{4}} e^{\tau \mu_j} &\gtrsim \|M_{\phi_j}\| \simeq \|U(i \tau)^{-1} \circ M_{\phi_j} \circ U(i \tau)\| \\ &\geq \left| \left\langle U(i \tau)^{-1} M_{\phi_j} U(i \tau) \phi_k, \phi_\ell \right\rangle \right| \\
    &= e^{- \tau \mu_k} \left| \left\langle U(i \tau)^{-1}  \mathcal A(\tau) (\phi_j \phi_k), \phi_\ell \right\rangle \right|\\
    &= e^{- \tau (\mu_k - \mu_\ell)} |c_{jk\ell}|,
\end{align*}
proving the claim.
\end{proof}

\subsection{Equivalent decay conditions} Conditions (c) and (d) in \Cref{thm:tech_real_analytic} may be obtained from condition (b) by the following straightforward arguments.

\begin{proof}[Proof of \Cref{thm:tech_real_analytic}] (b) $\Rightarrow$ (c).
 When $\mu_\ell>\mu_j+\mu_k+R$, condition (b) gives
 \[
    |c_{jk\ell}| \lesssim (1+\mu_\ell)^b e^{-a(\mu_\ell-\mu_j-\mu_k)}.
 \]
Since $1+\mu_\ell\leq (1+\mu_j+\mu_k)(1+\mu_\ell-\mu_j-\mu_k)$, we may rewrite this after multiplying by $(1+\mu_\ell)^n$ and absorbing a term into the exponential as
\[
    (1+\mu_\ell)^n|c_{jk\ell}|\lesssim (1+\mu_j+\mu_k)^{b+n}e^{-aR/2}.
\]
Rearranging, squaring, then summing over $\ell$ and applying Weyl's law gives
\[
    \sum_{\substack{\ell\textup{ with}\\ \mu_\ell>
            \mu_j+\mu_k+R}} |c_{jk\ell}|^2 \lesssim (1+\mu_j+\mu_k)^{2(b+n)}e^{-aR}\sum_{\ell=0}^\infty (1+\mu_\ell)^{-2n}\lesssim (1+\mu_j+\mu_k)^{2(b+n)}e^{-aR},
\]
and the result follows.

\medskip \noindent (c) $\Rightarrow$ (d). Take $a,b, C$ to be the constants coming from condition (c). We split the proof into two cases. First, when $|\mu_k-\mu_\ell|\leq 2\mu_j$, we simply use the fact that
\[
    |c_{jk\ell}|\leq \lVert\phi_j \rVert_\infty\leq e^{a\mu_j}\lVert\phi_j \rVert_\infty e^{-a|\mu_k-\mu_\ell|/2}.
\]
\noindent When instead $|\mu_k-\mu_\ell|> 2\mu_j$, set $R=|\mu_k-\mu_\ell|/2$. If $\mu_\ell>\mu_k$, then
\[
    \mu_\ell>\mu_j+\mu_k+R
\]
so, by (c),
\begin{align*}
    |c_{jk\ell}|\leq  \left( \sum_{\substack{r\textup{ with}\\ \mu_r>
            \mu_j+\mu_k+R}} |c_{jkr}|^2 \right)^{1/2}&\leq C(1+\mu_j+\mu_k)^b e^{-aR}\\
            &\leq C(1+\mu_j+\mu_k)^b e^{-a|\mu_k-\mu_\ell|/2}.
\end{align*}
If $\mu_k>\mu_\ell$, we note first that $\overline{c_{jk\ell}}=\langle\overline{\phi_j}\phi_\ell, \phi_k\rangle$. Expanding $\overline{\phi_j}$ in our basis and repeating the argument above for each term in this expansion, which necessarily has finitely many terms, yields the result.
\end{proof}

\subsection{Recovering real-analyticity from decay} It remains to show that condition (d) implies that the metric is real-analytic. We record the following lemma, which completes the proof, separately, as it shall also be employed in the next section.

\begin{lemma}\label[lemma]{lem:real_analytic_charts}
    Suppose that for every eigenfunction $\phi$, there exist constants $r, M>0$ such that, for every unit-speed geodesic $\gamma$, the function $\phi\circ \gamma$ extends holomorphically to the strip $\{z\in \mathbb{C}\,:|\Im(z)|<r\}$ with absolute value bounded by $M$. Then, the metric $g$ is real-analytic and all eigenfunctions are real-analytic.
\end{lemma}

\begin{proof}
    Fix a point $p\in M$ and an eigenfunction $\phi$, and define $\Phi=\phi\circ \exp_p\colon T_pM\to \mathbb{C}$. Standard Cauchy estimates give
    \[
        |D^k\Phi(0)(v,\dots,v)|\leq Mk!\, (|v|/r)^k
    \]
    for $v\in T_pM\cong \mathbb{R}^n$, so, by the polarization identity,
    \[
        \lVert D^k\Phi(0) \rVert\leq \frac{k^k}{k!}\sup_{|v|\leq 1}|D^k\Phi(0)(v,\dots,v)|\leq Mk! (e/r)^k.
    \]
    It follows that the Taylor series of $\Phi$ converges, hence $\Phi$ is real-analytic, in a neighborhood of $0$. Now, coordinate charts whose components are eigenfunctions cover the manifold, and it thus follows from our assumption that these charts are real-analytic with real-analytic transition maps.

    To see that $g$ is real-analytic when written in these coordinates, fix a chart with coordinate $x$ and write
    \[
        \Delta=\sum_{|\alpha|\leq 2}a_\alpha(x)\partial^\alpha.
    \]
    Then, with $k$ sufficiently large, choose eigenfunctions $\phi_1,\dots,\phi_k$ with eigenvalues $\lambda_1,\dots, \lambda_k$ so that we may solve for the coefficients $a_\alpha(x)$ from the system
    \[
        \sum_{|\alpha|\leq 2}a_\alpha(x)(\partial^\alpha\phi_j(x))=\lambda_j\phi_j(x), \quad 1\leq j\leq k.
    \]
    This shows the coefficients of the Laplacian are real-analytic. The metric coefficients may be recovered from the principal symbol of the Laplacian, and the result follows.
\end{proof}

\begin{proof}[Proof of \Cref{thm:tech_real_analytic}] (d) $\Rightarrow$ (a).
    This follows immediately from combining \Cref{prop:decay-to-geodesic} with \Cref{lem:real_analytic_charts}.
\end{proof}

\section{Characterizing flatness}\label{sec:flatness}

We turn to the proof of \Cref{thm:flat}. To begin, we establish the following fact, which is analogous to \Cref{lem:real_analytic_charts} and will allow us to restrict our attention to local considerations.

\begin{lemma}\label[lemma]{lem:finite_exponential_charts}
    Suppose that for every eigenfunction $\phi$, there exists a constant $N>0$ such that, for every geodesic $\gamma$, the function $\phi\circ \gamma$ is a linear combination of at most $N$ complex exponentials. Then, the metric $g$ is real-analytic and all eigenfunctions are real-analytic.
\end{lemma}
\begin{proof}
    Fix a point $p\in M$ and an eigenfunction $\phi$, and define $\Phi=\phi\circ \exp_p\colon T_pM\to \mathbb{C}$. Choose $v_0\in T_pM$ such that $\Phi(tv_0)$, when written as a linear combination of complex exponentials, has a maximal number $m$ of terms. This property is preserved near $v_0$, so, for $\delta>0$ sufficiently small and $v\in \overline{B_\delta(v_0)}$, we may write
    \[
        \Phi(tv)=\sum_{j=1}^ma_j(v)e^{i\tau_j(v)t}
    \]
    with $a_j,\tau_j$ uniformly bounded by some constant $C$. It follows that
    \[
        |D^k\Phi(0)(v,\dots, v)|\leq  mC^{k+1}
    \]
    in this neighborhood, so, by the polarization identity centered at $v_0$, we have
    \[
        \lVert D^k\Phi(0) \rVert\leq \frac{k^k}{\delta^k k!}\sup_{v\in \overline{B_\delta(v_0)}}|D^k\Phi(0)(v,\dots,v)|\leq mC (eC/\delta)^k.
    \]
    The Taylor series of $\Phi$ therefore converges, hence $\Phi$ is real-analytic, in a neighborhood of $0$. The remainder of the proof proceeds exactly as that of \Cref{lem:real_analytic_charts}.
\end{proof}

\subsection{A torus action from sparsity} We first work under a ``one-sided'', weakened version of the condition from \Cref{thm:flat}, assuming only that for each $j \in \mathbb N$, there exists $N_j \in \mathbb N$ such that $\#\{\ell: c_{jk\ell}\neq0\} \leq N_j$ for every $k \in \mathbb N$. This has the remarkable consequence that the topological closure of each geodesic orbit in $S^\ast M$ is a torus. In fact, a dense open subset of $S^\ast M$ admits a real-analytic action by the group $\mathbb T^m = (\mathbb R/2\pi \mathbb Z)^m$.

\begin{proposition}\label[proposition]{prop:torus-action}
    If there exists an orthonormal eigenbasis $(\phi_j)_{j=0}^\infty$ and integers $(N_j)_{j=0}^\infty$ such that every $\phi_j$ satisfies
    \[
        \#\{\mu_\ell \,:\, \langle \phi_j \phi_k,\phi_\ell \rangle \neq 0 \}\leq N_j
    \]
    for every $k$, then the topological closure of every geodesic orbit in $S^\ast M$ is a torus.

    Moreover, there exists a dense open subset $U \subseteq S^\ast M$, a real-analytic action $\rho: \mathbb T^m \times U \to U$ of the torus $\mathbb T^m$, $1 \leq m \leq n$ and a $\mathbb T^m$-invariant real-analytic map $\beta: U \to \mathbb R^m$ such that for every $\xi \in U$, $G_t \xi = \rho( t \beta(\xi))\, \xi$.

    For each element $\phi_j$ of the eigenbasis $(\phi_j)_{j=0}^\infty$ and all $\xi \in U$, the function $\tau \mapsto (\phi_j \circ \pi) (\rho(\tau)\, \xi)$ is a trigonometric polynomial in $\tau$ with at most $N_j$ terms.
\end{proposition}

\begin{proof}
    For simplicity, in this proof we assume without loss of generality that the eigenbasis $(\phi_j)_{j=0}^\infty$ under consideration is real valued.
    As a consequence of the local Weyl law, one may choose eigenfunctions $\phi_1,\ldots,\phi_L$ such that the map $\Phi: M \to \mathbb R^L$, $\Phi(p) = (\phi_1(p),\ldots,\phi_L(p))$ is an embedding. Write $\psi_j = \phi_j \circ \pi$, where $\pi: S^\ast M \to M$ is the canonical projection. By the local Weyl law, we may choose $\phi_1,\ldots,\phi_L$ so that $\nabla\phi_1, \ldots, \nabla\phi_L$ span $T_xM$ at all $x \in M$. Let $X$ denote the infinitesimal generator of the geodesic flow. It then follows that $S^\ast M$ embeds into $\mathbb R^{2L}$ via the map 
    \[\xi \mapsto (\psi_1(\xi), X \psi_1(\xi), \ldots, \psi_L(\xi), X \psi_L(\xi)).\]
    
    By \Cref{prop:sparse-to-geodesic}, $(\psi_j \circ G_t)(\xi)$ is a real-valued trigonometric polynomial in $t$ for every $\xi \in S^\ast M$.
    For each $j = 1,\ldots,L$, let $N_j$ be the largest number such that there exists $\xi_0 \in S^\ast M$ for which $(\psi_j \circ G_t)(\xi_0), (X\psi_j \circ G_t)(\xi_0), \ldots,(X^{N_j-1}\psi_j \circ G_t)(\xi_0)$ are linearly independent functions of $t$, but 
    \begin{align*}
        X^{N_j} \psi_j \circ G_t (\xi_0) = \sum_{k=0}^{N_j-1} a^j_k \, X^k (\psi_j \circ G_t) (\xi_0)
    \end{align*}
    for some (necessarily unique) $a^j_0,a^j_1,\ldots,a^j_{N_j-1} \in \mathbb R$. By \Cref{prop:sparse-to-geodesic} and \Cref{lem:finite_exponential_charts}, $M$, $g$ and hence all eigenfunctions and the vector field $X$ on $S^\ast M$ are real-analytic. Thus, $(\psi_j \circ G_t)(\xi),X (\psi_j \circ G_t)(\xi), \ldots,X^{N_j-1}(\psi_j \circ G_t)(\xi)$ remain linearly independent for $\xi$ in a dense open set $U \subseteq S^\ast M$ containing $\xi_0$. By maximality of $N_j$, for each $\xi \in U$ there exist unique $a^j_0(\xi),a^j_1(\xi),\ldots,a^j_{N_j-1}(\xi) \in \mathbb R$ such that 
    \begin{align*}
        X^{N_j}(\psi_j \circ G_t)(\xi) = \sum_{k=0}^{N_j-1} a^j_k(\xi) \, X^k (\psi_j \circ G_t) (\xi).
    \end{align*}
    From uniqueness of $(a^j_k(\xi))_{k=0}^{N_j-1}$ it follows that $a^j_k(\xi)$ is $G_t$-invariant, allowing us to choose $G_t$-invariant $U$. 
    This will remain true for the entire remaining proof: whenever $U$ is restricted to a possibly smaller open set, this can be done so $U$ remains dense and $G_t$-invariant.

    It is now easy to check that $(a^j_k(\xi))_{k=0}^{N_j-1}$ satisfies the above if and only if $(a^j_k(\xi))_{k=0}^{N_j-1}$ is the unique solution to 
    \begin{align*}
        \begin{pmatrix}
            \psi_j & X \psi_j & \dots & X^{N_j-1} \psi_j \\
            X\psi_j & X^2 \psi_j & \dots & X^{N_j} \psi_j \\
            \vdots & \vdots & \ddots & \vdots \\
            X^{N_j-1}\psi_j & X^{N_j} \psi_j & \dots & X^{2N_j-2} \psi_j \\
        \end{pmatrix} \begin{pmatrix}
            a^j_0 \\ a^j_1 \\ \vdots \\ a^j_{N_j-1}
        \end{pmatrix} = \begin{pmatrix}
            X^{N_j} \psi_j \\ X^{N_j+1} \psi_j \\ \vdots \\ X^{2N_j-1} \psi_j
        \end{pmatrix},
    \end{align*}
    which, by Cramer's rule, inherits the regularity of $M$, $X$ and $\psi_j$. Thus, by \Cref{sec:real-analyticity}, $a^j_k(\xi)$ is analytic in $\xi$ on the open set $U$.

    Let $N = \max(N_1,\ldots,N_L)$. Consider the embedding $\Psi: S^\ast M \to \mathbb R^{LN}$, 
    \begin{align*}
        \Psi(\xi) = \left( \psi_1(\xi),X\psi_1(\xi),\ldots,X^{N-1} \psi_1(\xi),\ldots, \psi_L(\xi),X \psi_L(\xi),\ldots,X^{N-1} \psi_L(\xi)\right).
    \end{align*}
    By the computations above, for fixed $\xi$, we have
    \begin{align*}
        \Psi_\ast X|_{\xi} = \mathbf B(\xi)\Psi(\xi),\qquad
        \mathbf B(\xi) = \operatorname{diag}(B_1(\xi),\ldots,B_L(\xi)),
    \end{align*}
    with
    \begin{align*}
        B_j(\xi) = \begin{pmatrix}
            0 & 1 & \cdots & 0 \\
            0 & 0 & \cdots & 0 \\
            \vdots & \vdots & \ddots & \vdots \\
            0 & 0 & \cdots & 1 \\
            b^j_0(\xi) & b^j_1(\xi) & \cdots & b^j_{N-1}(\xi)
        \end{pmatrix},
    \end{align*}
    where, to account for the different values of $N_1,\ldots,N_L$, we set $b^j_k(\xi) = a^j_{k-N+N_j}(\xi)$ if $k \geq N-N_j$ and $b^j_k(\xi) = 0$ otherwise. Since $b^j_k(\xi)$ is constant along orbits of $X$, it follows that for every $\xi \in U$,
    \begin{align*}
        \Psi(G_t(\xi)) = \exp(t \mathbf B(\xi)) \Psi(\xi).
    \end{align*}
    Since $\Psi(G_t(\xi))$ stays bounded for all $t \in \mathbb R$, it follows that there exist purely imaginary eigenvalues $i \alpha_1(\xi),\ldots,i \alpha_\ell(\xi) \in \spec(\mathbf B(\xi)) \cap i \mathbb R$ and corresponding eigenvectors $V_1(\xi),\ldots,V_\ell(\xi)$ such that 
    \begin{align*}
        \Psi(G_t(\xi)) = e^{i \alpha_1(\xi) t} V_1(\xi) + \ldots + e^{i \alpha_\ell(\xi) t} V_\ell(\xi).
    \end{align*}
    If $\alpha_1,\ldots,\alpha_\ell$ were linearly independent over $\mathbb Q$, Kronecker's theorem on simultaneous Diophantine approximation would imply
    \begin{align*}
        \overline{\{ \Psi(G_t(\xi)): t\in \mathbb R \}} = \{ e^{i t_1} V_1(\xi) + \ldots + e^{i t_\ell} V_\ell(\xi) : (t_1,\ldots,t_\ell) \in \mathbb T^\ell \}
    \end{align*}
    for all $\xi \in U$, after possibly shrinking $U$. In the present setting, that need not be the case. At the very least, if $(\phi_j)_{j=0}^\infty$ consists of real valued eigenfunctions, $\Psi(G_t(\xi))$ is real, and for every $\alpha_j > 0$ there must exist $\alpha_{j'}$ with $\alpha_{j'} + \alpha_j = 0$. Further linear dependencies over $\mathbb Q$ may exist, and in the case of flat manifolds, indeed occur for $\xi$ in a dense subset of $S^\ast M$.
    
    We remedy this situation as follows. Since $\mathbf B(\xi)$ is real-analytic, on a dense open subset of $U$ the frequencies $\alpha_1,\ldots,\alpha_\ell$ are real-analytic functions of $\xi$ as well, and any linear relation between them holds either everywhere or on a nowhere dense closed subset of $U$. By Baire's theorem, there exists a comeagre subset $B$ of $U$ such that any $\mathbb Q$-linear relation between $\alpha_1,\ldots,\alpha_\ell$ holds identically on $B$. Eliminating these dependencies yields real numbers $\beta_1(\xi),\ldots,\beta_m(\xi) \in \mathbb R$, which are linearly independent over $\mathbb Q$, and $(\kappa_{j,k})_{j,k} \in \mathbb Z^{\ell \times m}$ such that $\alpha_j(\xi) = \sum_{k} \kappa_{j,k} \beta_k(\xi)$. It then follows that 
    \begin{align*}
        \overline{\{ \Psi(G_t(\xi)): t\in \mathbb R \}} = \Bigg\{ \sum_{j=1}^\ell e^{i \kappa_{j,1} t_1 + \ldots + i \kappa_{j,m} t_m} V_j(\xi) : (t_1,\ldots,t_m) \in \mathbb T^m \Bigg\}
    \end{align*}
    for all $\xi$ in $B$. Denote $T_\xi := \overline{\{ \Psi(G_t(\xi)): t\in \mathbb R \}}$.
    
    Since the tori $T_\xi$ arise as minimal sets of an autonomous flow, for any $\xi,\xi' \in U$ we have the dichotomy $T_{\xi} \cap T_{\xi'} \in \{\varnothing, T_{\xi}\}$. 
    It is not immediately clear that $T_\xi \cap U = T_\xi$. However, this does follow from the fact that the flow of $X$ on $T_\xi$ is conjugate to the flow of an irrational vector field on $\mathbb T^m$, and $U$ is invariant under the flow of $X$. Hence, $U = \bigcup_{\xi \in U} T_\xi$.
    
     Since the numbers $\alpha_j$ are unique and constant on orbits $G_t(\xi)$, and these orbits are dense in $T_\xi$, we infer that $\alpha_j$ is constant on $T_\xi$. The vector fields
    \begin{align*}
        Y_k(\xi) := i \sum_{j = 1}^\ell \kappa_{j,k} V_j(\xi),
    \end{align*}
    are tangential to $T_\xi$ wherever $\xi \in B$, and by construction commute with $X$ there. Since $B$ is dense, the action $\rho: \mathbb T^m \times U \to U$ generated by $Y_1,\ldots,Y_m$ commutes with the geodesic flow $X$ on $U$. Also by construction, 
    \begin{align*}
        G_t(\xi) = \rho(t \beta_1(\xi), \ldots, t \beta_m(\xi)) \, \xi.
    \end{align*}

    The construction of this torus action imposes no restrictions on the dimension $m$. However, the following argument shows that $m \leq n$. We begin with the well-known observation that the tautological form on $T^\ast M$ restricts to a contact form $\lambda$ on $S^\ast M$ -- that is, a form with the property that $d \lambda$ is a nondegenerate two-form on $\ker \lambda$ -- which is furthermore invariant under the flow of $X$ \cite[Section 10.5.2]{Zelditch2017}. Now pick $\xi \in B$, so that $T_\xi$ is $m$-dimensional and write $\rho_\xi: \mathbb T^m \to U$, $\rho_\xi(\tau) = \rho(\tau) \, \xi$. Invariance under $X$ translates to the fact that $\rho_\xi^\ast \lambda$ has constant coefficients in the standard coordinate system on $\mathbb T^m$, since the pull-back of $X$ under $\rho_\xi$, being an irrational vector field on $\mathbb T^m$, generates a dense subset of translations on $\mathbb T^m$. Combining $\rho_\xi^\ast d\lambda = 0$ with the fact that $d\lambda$ is a nondegenerate two-form on $\ker \lambda$, we infer that $\ker \lambda \cap d\rho_\xi (T \, \mathbb T^m)$ has dimension at most $n-1$, and hence, that $m \leq n$.
    
    To obtain the last statement, it suffices to invoke the following well-known fact from the theory of almost periodic functions: A continuous function $f: \mathbb T^m \to \mathbb C$ such that $f(t \beta_1,\ldots,t\beta_m)$ is a trigonometric polynomial in $t$ with at most $N$ terms for some $\mathbb Q$-linearly independent $\beta_1,\ldots,\beta_m$ is a trigonometric polynomial on $\mathbb T^m$ with at most $N$ terms \cite[Chapter II, Section 8]{BrockerTomDieck1985}.
\end{proof}

\subsection{Consequences of uniform sparsity} We now complete the proof of  \Cref{thm:flat}. For the remainder of \Cref{sec:flatness}, we assume that there exists an orthonormal eigenbasis $(\phi_j)_{j=0}^\infty$ and $N \in \mathbb N$ such that for every $j,k \in \mathbb N$,
    \[
        \#\{\mu_\ell \,:\, \langle \phi_j \phi_k,\phi_\ell \rangle \neq 0 \}\leq N.
    \]
    This allows us to deduce our second main result, namely that $(M,g)$ is flat and hence a quotient of Euclidean $n$-space by a cocompact subgroup of its group of isometries.
    
    We begin by considering the action, constructed in \Cref{prop:torus-action}, of $\mathbb T^m$ on an open subset $U \subset S^\ast M$. We will show that under the additional assumption that $N_j \leq N$, we may infer that $m=n$.

    \begin{proposition}\label[proposition]{prop:torus-dimension}
    Let $(\phi_j)_{j=0}^\infty$, $(N_j)_{j=0}^\infty$, $U$ and $\rho: \mathbb T^m \times U \to U$ be as in \Cref{prop:torus-action}. If there exists $N \in \mathbb N$ such that $N_j \leq N$ for all $\phi_j$ in the sequence $(\phi_j)_{j=0}^\infty$, then $m=n$, and the map $\rho_\xi: \mathbb T^m \to M$, $\rho_\xi(\tau) = \pi(\rho(\tau) \xi)$ is a covering map for all $\xi \in U$.
\end{proposition}
    
\begin{proof}
    Consider the infinitesimal generators $Y_1,\ldots,Y_m \in \Gamma(T U)$ of the torus action, that is,
    \begin{align*}
        Y_\ell = \tfrac{d}{dt}\big|_{t=0} \, \rho(t e_\ell), \ \ell = 1,\ldots,m.
    \end{align*}
    The bulk of this proof consists of showing that, locally near any $\xi_0 \in U$, $\pi_\ast Y_\ell|_\xi$ depends only on the basepoint $\pi(\xi)$. This then immediately implies $m \geq n$: Since $X = \sum_{k=1}^m \beta_k Y_k$, 
    \begin{align*}
        \xi^\# = \pi_\ast X|_{\xi} \in \mathrm{span}(\pi_\ast Y_1|_{\xi},\ldots,\pi_\ast Y_{m}|_{\xi}),
    \end{align*}
    so the vectors $\{ \pi_\ast Y_1|_{\xi},\ldots,\pi_\ast Y_{m}|_{\xi} : \xi \in S_p^\ast M\}$ span $T_p M$ whenever $p \in \pi(U)$. Since $m \leq n$ by the contact-geometric argument from \Cref{prop:torus-action}, $m=n$ as claimed. It furthermore follows that $\tau \mapsto \pi(\rho(\tau) \xi)$ is a submersion between equidimensional closed manifolds, i.e., a covering map. For the remainder of the proof, fix one generator $Y_\ell$.
    
    Let $(\phi_{j_k})_{k=0}^\infty$ be a subsequence giving rise to a semiclassical defect measure $\mu$ on $S^\ast M$. We will write $h_k = 1/\mu_{j_k}$ and treat $(\phi_{j_k})_{k=0}^\infty$ as an $h$-tempered semiclassical family of distributions $\phi_h$, usually dropping the index $k$ altogether. 
    
    To understand the semiclassical defect measure on $T^\ast (S^\ast M)$ associated to the family $\phi_h \circ \pi$, we first observe that for every $x \in M$, the fibers of the conormal bundle $N^\ast (S^\ast_x M)$ are in natural bijection with $T^\ast_x M$ by pull-back. In other words, given $\xi \in S^\ast_x M$, we can identify $T^\ast_x M$ with $N^\ast_\xi (S^\ast_x M)$ by associating $\alpha \in T^\ast_x M$ with $(d\pi|_\xi)^\ast \alpha \in T^\ast_\xi( S^\ast M)$ -- note that $(d\pi|_\xi)^\ast \alpha$ vanishes on $T_\xi(S^\ast_x M)$, since $d\pi_\xi ( T_\xi(S^\ast_x M)) = 0$, and hence $(d\pi|_\xi)^\ast \alpha \in N_\xi^\ast (S^\ast_x M)$. If $\beta \in N^\ast_\xi S^\ast_x M$ is given by $\beta = (d\pi|_\xi)^\ast \alpha$, we also write $\pi_\ast \beta = \alpha$. For fixed $\alpha \in T^\ast M$, the pre-image $\pi_\ast^{-1}(\alpha)$ can be identified with $S^\ast_{\pi(\alpha)} M$ and equipped with the uniform measure $\sigma$ on $S^\ast_{\pi(\alpha)} M$. One can then check that the family $\phi_h \circ \pi$ gives rise to the semiclassical defect measure $\nu$ implicitly defined by 
    \begin{align*}
        \int_{T^\ast S^\ast M} f \, d\nu &= \int_{\alpha \in T^\ast M} \left( \int_{\pi_\ast^{-1}(\alpha)} f d\sigma \right) d\mu(\alpha) .
    \end{align*}
    One may view $\nu$ as the pull-back measure along the fibration $\pi_\ast: N^\ast S^\ast M \to T^\ast M$. As an important consequence,
    \begin{align}
    \label{eq:support_after_lift}
        \supp(\nu) = \left\{ \beta \in N^\ast S^\ast M: \beta = (d \pi|_\xi)^\ast \alpha \text{ for } \alpha \in \supp \mu, \xi \in S^\ast M \right\} = \pi_\ast^{-1}(\supp(\mu)).
    \end{align}
    We will not use any property of $\nu$ other than this restriction on its support.
    
    Next, we use the torus action generated by $Y_1,\ldots,Y_m$ to obtain another strong restriction on $\supp(\nu)$. Fix an eigenfunction $\phi$ from  the sequence $\phi_h$. For each $\xi \in U$, consider the function $\psi_\xi \in C^\infty(\mathbb T^m)$ given by $\psi_\xi(\tau) = (\phi \circ \pi)(\rho(\tau)\xi)$. For fixed $\xi$, \Cref{prop:torus-action} implies that $\widehat{\psi_\xi}$ is supported on a finite collection $K_\xi \subseteq \mathbb Z^m$ with $|K_\xi| \leq N$. After restricting $U$ to the dense open subset where $|K_\xi|$ attains a local maximum, $K_\xi$ must be constant, since $\widehat{\psi_\xi}(\kappa)$ depends real-analytically on the covector $\xi$. Choose a polynomial $p \in \mathbb R[y]$ of degree at most $N$ such that $p(h \, \kappa \cdot e_\ell) = 0$ for all $\kappa \in K$, normalized such that $\max_{j} |p_j| = 1$. It follows that $p(-i h \partial_{\ell}) \psi_\xi = 0$ for all $\xi$ sufficiently close to $\xi_0$, and hence $p(- i h Y_\ell) (\phi \circ \pi) = 0$ on a neighborhood of $\xi_0$ in $S^\ast M$. By analyticity, $p(- i h Y_\ell) (\phi \circ \pi) = 0$ on the entire set $U\subseteq S^\ast M$ on which $Y_\ell$ is defined. Repeat this procedure for each function $\phi_h$ and denote the resulting polynomial by $p_h$.

    After passing to a subsequence, we may assume that $p_h \to p \in \mathbb R[y]$.  Since $\max_j |p_j| = 1$, $p$ does not vanish identically. Let $P_h = (1 + h^2 \Delta_{S^\ast M})^{-\frac N2} \circ p(-i h Y_\ell)$. Then
    \begin{align*}
        \| P_h (\phi_h \circ \pi) &\|_{L^2(S^\ast M)} = \| p(-i h Y_\ell) (\phi_h \circ \pi)\|_{H_h^{-N}(S^\ast M)} \\ &\lesssim \left\| \big( p(-i h Y_\ell)-p_h(-i h Y_\ell)\big) (\phi_h \circ \pi) \right\|_{H_h^{-N}(S^\ast M)} \!+\! \| p_h(-i h Y_\ell) (\phi_h \circ \pi) \|_{H_h^{-N}(S^\ast M)} \\
        &\lesssim \sum_{j=0}^N \left|p_j - (p_h)_j\right| \|(-i h Y_\ell)^j (\phi_h \circ \pi) \|_{H_h^{-N}(S^\ast M)} \\
        &\lesssim \Bigg( \sum_{j=0}^N |p_j - (p_h)_j| \Bigg) \|  (\phi_h \circ \pi) \|_{L^{2}(S^\ast M)} = o(1).
    \end{align*}
    By \cite[Theorem 5.3]{Zworski2012}, it follows that $\supp(\nu) \subseteq \{\zeta \in T^\ast S^\ast M: p(\left\langle \zeta,  Y_\ell\right\rangle ) = 0\}$. From \eqref{eq:support_after_lift}, we know that $\supp(\nu) \subseteq N^\ast S^\ast M$, and for $\zeta \in N^\ast S^\ast M$, it holds that $\left\langle \zeta,  Y_\ell\right\rangle = \left\langle \pi_\ast \zeta, \pi_\ast Y_\ell\right\rangle$. Combining this observation with \eqref{eq:support_after_lift} yields
    \begin{align}
    \label{eq:support_after_lift2}
        \pi_\ast^{-1}(\supp(\mu)) = \supp(\nu) \subseteq \{\beta \in N^\ast S^\ast M: \left\langle \pi_\ast \beta, \pi_\ast Y_\ell\right\rangle \in p^{-1}(0)\}.
    \end{align}

    This immediately restricts the support of the original semiclassical defect measure $\mu$:
    \begin{align}
    \label{eq:support_after_lift3}
        \supp(\mu) \subseteq \pi_\ast(\pi_\ast^{-1}(\supp(\mu))) \subseteq \{\alpha \in T^\ast M: \left\langle \alpha, \pi_\ast Y_\ell\right\rangle \in p^{-1}(0)\}.
    \end{align}
    Fix a point $\alpha \in \supp(\mu)$, and vary $\xi \in S^\ast_{\pi(\alpha)} M$. Since $\pi_\ast Y_\ell|_\xi$ depends continuously on $\xi$ and $p^{-1}(0)$ is discrete, $\left\langle \alpha, \pi_\ast Y_\ell|_{\xi}\right\rangle$ is a constant independent of $\xi$.

    Recall the well known fact (\Cref{lem:semiclassical-measures-are-dense}) that the union of $\supp(\mu)$ over all semiclassical measures arising from subsequences of $(\phi_j)_{j=0}^\infty$ is dense in $S^\ast M$. If $\alpha \in T^\ast M$ is contained in the support of some semiclassical measure $\mu_\alpha$, then the calculations above show that $\left\langle \alpha, \pi_\ast Y_\ell|_{\xi} \right\rangle = c_\alpha$ for some constant $c_\alpha \in \mathbb R$ independent of $\xi \in S^\ast_{\pi(\alpha)} M$. By continuity, the function 
    \begin{align*}
        (\alpha,\xi,\xi' ) \mapsto \left\langle \alpha, \pi_\ast Y_\ell|_{\xi} - \pi_\ast Y_\ell|_{\xi'} \right\rangle
    \end{align*}
    vanishes identically on $T^\ast M \times S^\ast M \times S^\ast M$, that is, $\pi_\ast Y_\ell$ does not depend on $\xi$.
\end{proof}

\begin{proof}[Proof of \Cref{thm:tech_flat}]
(a) $\Rightarrow$ (b). This follows from \cite[Lemma 4.1]{SchildkrautSpeciel2026}.

\medskip \noindent (b) $\Rightarrow$ (c). This is clear.

\medskip \noindent (c) $\Rightarrow$ (a). Consider a subsequence $\phi_h$ of $(\phi_j)_{j=0}^\infty$ giving rise to a semiclassical measure $\mu$ on $T^\ast M$. Fix $\xi \in U$ and consider $\rho_\xi: \mathbb T^m \to M$, $\rho_\xi(\tau) = \pi(\rho(\tau)\, \xi)$ as introduced in \Cref{prop:torus-dimension}. Since each $\phi_h \circ \rho_\xi$ is a trigonometric polynomial with at most $N$ terms on $\mathbb T^m$, the support of the semiclassical measure $\nu$ on $T^\ast \mathbb T^m$ induced by $\phi_h \circ \rho_\xi$ is invariant under translation. On the other hand, $\supp(\nu) = d\rho_\xi^\ast(\supp(\mu)) \subseteq d\rho_\xi^\ast(S^\ast M)$. Since
\begin{align*}
    S^\ast M = \overline{\bigcup_{\mu} \supp(\mu)},
\end{align*}
where the union ranges over all semiclassical defect measures $\mu$ arising from $(\phi_j)_{j=1}^\infty$, it follows that
\begin{align*}
    d\rho_\xi^\ast(S^\ast M) = \overline{\bigcup_{\mu} \supp(\nu)},
\end{align*}
and hence that $d\rho_\xi^\ast(S^\ast M)$ is translation-invariant. It follows that the metric $\rho_\xi^\ast g$ on $\mathbb T^m$ is translation-invariant, hence $M$ is flat.

It remains to show that, when $(M,g)$ is flat, the minimal $N$ attained in conditions (b) and (c) is exactly the order of its holonomy group, which we denote by $H$. That the minimal such $N$ is bounded above by $|H|$ is shown in the proof of \cite[Lemma 4.1]{SchildkrautSpeciel2026}. The reverse inequality follows from the proof of \cite[Lemma 4.2]{SchildkrautSpeciel2026}, after replacing (in the notation of \cite{SchildkrautSpeciel2026}) $H_1$, the group generated by $H$ and $-\operatorname{Id}$, and
\[
    \Re(u_{\xi_0}), \quad\textrm{where}\quad u_\xi(x)=\frac{1}{|H|}\sum_{h\in H}e^{2\pi i \langle \xi, h(x)\rangle}
\]
and $\xi_0$ is a lattice element whose $H$ orbit has maximal size, by $H$ and $u_{\xi_0}$ instead.
\end{proof}

\section{Non-locally homogeneous metrics on \texorpdfstring{$\mathbb S^7$}{S7}}\label{sec:weird-spheres}

We construct a family of metrics on $\mathbb{S}^7$ which satisfies the properties described in \Cref{thm:tech_nonhomogeneous}. To this end, consider $\mathbb{S}^7$ as a subset of $\mathbb{R}^8\cong \mathbb{H}^2$, where $\mathbb{H}$ denotes the space of quaternions with orthonormal basis $1, q_1, q_2$ and $q_3$. Recall that $\mathbb{H}$ carries a noncommutative multiplicative structure defined by $q_j^2=-1$ and $q_1q_2=q_3$, a conjugation defined by $\overline {1}=1$ and $\overline {q_j}=-q_j$, and an inner product defined by $\langle u,v\rangle _\mathbb{H}=\Re(u\overline v)$. We may therefore identify $\mathbb{S}^7\cong \{(u,v)\in \mathbb{H}^2\,\colon\, |u|^2+|v|^2=1\}$, and we fix $g_1$ to be the standard metric on $\mathbb{S}^7$.

Define for $j=1,2,3$ the linear map
\[
    M_j\colon \mathbb{H}^2\to \mathbb{H}^2\quad\textrm{by}\quad M_j(u,v)=(q_ju,-q_jv),
\]
and observe $\langle M_j(u,v),(u,v)\rangle =0$, hence the vector fields given by $X_j(u,v)=M_j(u,v)$, after identifying the tangent space with the origin by translation, are tangent to $\mathbb{S}^7$. These vector fields generate the rotations $e^{tM_j}$ and are thus Killing fields of $g_1$. Furthermore, a simple computation confirms that they are orthonormal at every point.

Denote by $V$ the subbundle of $T\mathbb{S}^7$ spanned by $X_1, X_2, X_3$, and by $P_V$, $P_{V^\perp}$ the orthogonal projections onto $V$ and $V^\perp$, respectively. Introduce for $t>0$ the metrics $g_t$ defined by
\begin{equation}\label{eq:metric_def}
    g_t(Y,Z)=t^{4/7}g_1(P_VY, P_VZ)+t^{-3/7}g_1(P_{V^\perp}Y, P_{V^\perp}Z)
\end{equation}
for tangent vectors $Y,Z$. We shall prove that, when $t$ is sufficiently close to $1$, this family has the desired properties.

\begin{proof}[Proof of \Cref{thm:tech_nonhomogeneous}]
    Since $\dim V=3$ and $\dim V^\perp=4$, we have
    \[
        d\textrm{vol}_{g_t}= (t^{4/7})^{3/2}(t^{-3/7})^{4/2}d\textrm{vol}_{g_1}=d\textrm{vol}_{g_1}.
    \]
    Furthermore, the metrics $g_t$ are themselves real-analytic, and setting $t=1$ recovers $g_1$. It remains to prove claims (c) and (d) of the theorem.

    To establish claim (c), write $\Delta_t$ for the Laplacian corresponding to $g_t$ and define
    \[
        \Delta_V=-\sum_{j=1}^3X_j^2.
    \]
    Since the vector fields $X_j$ are Killing fields, the operators $\Delta_1$ and $\Delta_V$ commute. It follows that $\Delta_V$ restricts to a well defined map on the eigenspaces of $\Delta_1$. Denote these eigenspaces by $E_d$, recalling that $E_d$ consists of restrictions of homogeneous harmonic polynomials of degree $d$ in $\mathbb{R}^8$ to $\mathbb{S}^7$. Now, since $\Delta_V$ is well defined and symmetric on $E_d$, we may construct an eigenbasis shared by $\Delta_1$ and $\Delta_V$. A simple computation reveals
    \begin{equation}
        \Delta_t=t^{-4/7}\Delta_V+t^{3/7}(\Delta_1-\Delta_V)=t^{3/7}(\Delta_1+(t^{-1}-1)\Delta_V),
    \end{equation}
    hence this eigenbasis in fact simultaneously diagonalizes all operators $\Delta_t$, and the claim follows.

    Finally, we prove claim (d). To begin, we show that, with $F(u,v)=|u|^2-|v|^2$, we have
    \begin{equation}\label{eq:intersection_of_espaces}
        E_2\cap \ker \Delta_V=\textrm{span}_\mathbb{C}F.
    \end{equation}
    Indeed, $F$ is the restriction of a homogeneous harmonic polynomial which is constant when restricted to $\{(u,0)\,\colon |u|\textrm{ is constant}\}$ or $\{(0,v)\,\colon |v|\textrm{ is constant}\}$. Therefore, $F$ is invariant under the flow of $X_j$, and the right-to-left inclusion follows. Conversely, note for $f\in E_2\cap \ker \Delta_V$ that
    \[
        0=\langle \Delta_Vf,f\rangle=\sum_{j=1}^3\|X_jf\|^2.
    \]
    It follows that $f$ is constant along the flow of $X_j$, and so invariant under $e^{tM_j}$. Since $e^{(\pi/2)M_j}=M_j$ and $M_1M_2M_3(u,v)=(-u,v)$, it follows that $f(u,v)=f(-u,v)$. However, as $f\in E_2$, its mixed quadratic terms must vanish so that $f(u,v)=f(u,0)+f(0,v)$, now viewing $f$ as a homogeneous harmonic polynomial on $\R^8$. But invariance under $e^{tM_j}$ implies that each of these terms is constant on spheres in $\mathbb{R}^4$, so we must have
    \[
        f(u,v)=\alpha|u|^2+\beta|v|^2\quad\textrm{for some }\alpha, \beta\in \mathbb{C}.
    \]
    \Cref{eq:intersection_of_espaces} follows since $f$ is harmonic in $\R^8$ so $\alpha=-\beta$.

    Next, note
    \[
        \Delta_tF=t^{3/7}(\Delta_1+(t^{-1}-1)\Delta_V)F=16t^{3/7}F.
    \]
    When $t\neq 1$, it follows that the function $F$ is, up to scaling, the only eigenfunction of $\Delta_t$ in $E_2$ with eigenvalue $16t^{3/7}$, since the equation above forces any such eigenfunction to be in the kernel of $\Delta_V$. Now, in a small neighborhood of $16$, there are exactly $\dim E_2$ eigenvalues of $\Delta_t$ for $t$ sufficiently close to $1$ \cite[Chapter IV, Theorem 3.16]{Kato1995}. It follows that $16t^{3/7}$ is a simple eigenvalue with eigenfunction $F$ when $t$ is sufficiently close, but not equal, to $1$.

    Any isometry of $g_t$ commutes with $\Delta_t$, hence maps eigenfunctions to eigenfunctions and preserves their eigenvalues. Therefore, when $t\neq1$ is sufficiently close to $1$, any such isometry must map $F$ to $\pm F$. Deduce the isometry group of $g_t$ cannot act transitively on $\mathbb{S}^7$. It is shown in \cite[Section 4, Corollary]{Nomizu1962} that any locally homogeneous, simply connected, complete, connected Riemannian manifold must be homogeneous, from which the result now follows.
\end{proof}

We note that the analogous construction does not work on $\mathbb{S}^3$, instead yielding Berger's sphere, which is homogeneous \cite[Example 9.81]{Besse1987}.

\bibliographystyle{plain}
\bibliography{references}

\end{document}